\documentclass[12pt]{article}
\usepackage{amssymb,latexsym}
\def\ZZ{{\mathbb Z}}
\def\QQ{{\mathbb Q}}
\def\RR{{\mathbb R}}

\newtheorem{formula}{}[section]
\newtheorem{proposition}[formula]{Proposition}
\newtheorem{definition}[formula]{\indent Definition}
\newtheorem{corollary}[formula]{\indent Corollary}
\newtheorem{remark}[formula]{\indent Remark}
\newtheorem{lemma}[formula]{\indent Lemma}
\newtheorem{theorem}[formula]{\indent Theorem}

\newtheorem{example}[formula]{Example}

\def\thrm{\begin{theorem}}
\def\thrml#1{\begin{theorem}\label{#1}}
\def\ethrm{\end{theorem}}
\def\rmrk{\begin{remark}}
\def\rmrkl#1{\begin{remark}\label{#1}}
\def\ermrk{\end{remark}}
\def\dfntn{\begin{definition}}
\def\dfntnl#1{\begin{definition}\label{#1}}
\def\edfntn{\end{definition}}
\def\nmrt{\begin{enumerate}}
\def\enmrt{\end{enumerate}}

\def\qtn{\begin{equation}}
\def\qtnl#1{\begin{equation}\label{#1}}
\def\eqtn{\end{equation}}
\def\lmm{\begin{lemma}}
\def\lmml#1{\begin{lemma}\label{#1}}
\def\elmm{\end{lemma}}
\def\crllr{\begin{corollary}}
\def\crllrl#1{\begin{corollary}\label{#1}}
\def\ecrllr{\end{corollary}}

\begin{document}
\title{}
\date{}
\maketitle
%\nopagenumbers
%\begin{titlepage}
%\title
\vspace{-0,1cm} \centerline{\bf Tropical recurrent sequences}
%\centerline{\bf AND RELATIVE KOLMOGOROV COMPLEXITY}
\vspace{7mm}
\author{
\centerline{Dima Grigoriev}
%\\[-1pt]
\vspace{3mm}
%\centerline{$^1$ St.Petersbourg University,
% Universitetskaya nab., 7/9,}
%\centerline{ St.Petersbourg,
% 199164, RUSSIA}
%\vspace{3mm}
%Dima Grigoriev \\[-1pt]
\centerline{CNRS, Math\'ematique, Universit\'e de Lille, Villeneuve
d'Ascq, 59655, France} \vspace{1mm} \centerline{e-mail:\
dmitry.grigoryev@math.univ-lille1.fr } \vspace{1mm}
\centerline{URL:\ http://en.wikipedia.org/wiki/Dima\_Grigoriev} }
%\date{}
%\maketitle

\begin{abstract}
Tropical recurrent sequences are introduced satisfying a given vector (being a
tropical counterpart of classical linear recurrent sequences). We consider the
case when Newton polygon of the vector has a single (bounded) edge. In this case there are periodic tropical recurrent sequences which are similar to classical
linear recurrent sequences. A question is studied when there exists a non-periodic tropical recurrent sequence satisfying a given vector, and partial answers
are provided to this question. Also an algorithm is designed which tests
existence of non-periodic tropical recurrent sequences satisfying a given
vector with integer coordinates. Finally, we introduce a tropical entropy of a vector, provide some bounds on it and extend this concept to tropical multivariable recurrent sequences.
\end{abstract}

{\bf keywords}: tropical recurrent sequence, periodic sequence, tropical entropy

\section*{Introduction}
A classical (linear) recurrent sequence $\{z_l\}_{l\in \ZZ}$ (e.~g. Fibonacci numbers) satisfies
conditions $\sum_{0\le i\le n} a_iz_{i+k}=0,\, k\in \ZZ,\, a_0\neq 0,\, a_n\neq 0$. It well known that the
linear space of all such sequences has dimension $n$ and can be explicitly
described via the roots of polynomial $\sum_{0\le i\le n} a_ix^i$ and the
derivatives in case of multiple roots.

We study {\it tropical recurrent sequences} satisfying similar tropical linear polynomials
$\min_{0\le i\le n}\{a_i+y_{i+k}\},\, k\in \ZZ$ where as it is adopted in tropical
algebra \cite{MS} we assume that the minimum is attained for at least two
different indices $0\le i_1(k) < i_2(k)\le n$ for each $k\in \ZZ$. We say that
in this case  a tropical recurrent sequence $y=\{y_i\in \RR \cup \{\infty\}\}_{i\in \ZZ}$ {\it satisfies vector $a:=(a_0,\dots,a_n)\in (\RR \cup \{\infty\})^{n+1}, a_0< \infty,\, a_n< \infty$}. One can treat
a tropical recurrent sequence as a solution of an (infinite) tropical
Macaulay matrix whose rows are obtained from vector $a$ by all possible shifts.  We mention that Macaulay matrix plays a key role in the tropical Nullstellensatz \cite{G}, \cite{GP}.

On the other hand, one can view tropical linear recurrence as a tropical analog of shift equations with constant coefficients (note that in the classical setting shift equations extend difference equations). Also one can obtain tropical recurrent sequences by means of tropicalizing linear recurrent sequences of the form $\sum _{0\le i\le n} A_iZ_{i+k}=0,\, k\in \ZZ$ over Puiseux series.

Throughout the paper (except for sections~\ref{five}, \ref{six}) we impose the requirement of {\it minimality} of tropical recurrent sequences: for any $j\in \ZZ$ one can not diminish $y_j$ keeping all the rest
$y_i,\, i\neq j$ without violation of being a tropical recurrent sequence
satisfying $a$.

A description of tropical recurrent sequences satisfying a given vector $a$ is more complicated than its classical counterpart, and we don't provide a
complete answer.

The collection of points $\{(i,a_i)\, : \, 0\le i\le n\} \subset \RR^2$ we call {\it Newton graph} of $a$.
A crucial feature of $a$ is its {\it Newton polygon} $P(a)$ on the plane which is the
convex hull of the vertical rays $\{(i,b\ge a_i),\, 0\le i\le n\}$. For each
(bounded)
edge of $P(a)$ with a slope $\sigma$ there exists a tropical recurrent sequence $\{y_j\}_{j\in \ZZ}$ satisfies $a$ where points $(j,y_j),\, j\in \ZZ$ are located on a line with
the slope $-\sigma$.
%When $P(a)$ has just a single (bounded) edge
There can be more
general periodic tropical recurrent sequences, so for some period $d\in \ZZ$
it holds $y_{j+d}-y_j=-\sigma d,\, j\in \ZZ$.

In a certain (informal) sense periodic tropical recurrent sequences are similar to classical recurrent sequences. On the other hand, for some vectors $a$ there exist
non-periodic tropical recurrent sequences satisfying $a$. We study for which
$a$ they exist.

In section~\ref{one}  vectors $a$ are considered such that all the finite points $(i,a_i),\, a_i<\infty$ from Newton graph lie  on a single (bounded) edge
of $P(a)$. We show that all the tropical recurrent sequences satisfying $a$ are
periodic iff the points $(i,a_i),\, a_i<\infty$ form an arithmetic progression
with some difference $d$. In the latter case any tropical recurrent sequence
satisfying $a$ has  period $d$.

In section~\ref{two}  a more general situation is studied when $P(a)$ has a
single (bounded) edge. First, we note that if points $(i,a_i)$ lying on this
edge do not form an arithmetic progression then there exists a non-periodic
tropical recurrent sequence satisfying $a$. When, on the contrary, these points form an
arithmetic progression with the difference 2 we prove that all the tropical
recurrent sequences satisfying $a$ are periodic iff points $(i,a_i)$ not lying
on the edge, also form an arithmetic progression.

In section~\ref{three} vectors $a$ are considered with an arithmetic progression having the difference 3 of points
$(i,a_i)$ lying on the (single) edge. We provide two examples of $a$: one having only periodic tropical recurrent sequences satisfying $a$, and another one
with non-periodic sequences. These two examples demonstrate that the existence
of non-periodic tropical recurrent sequences satisfying $a$ can not be expressed just in terms of arithmetic progressions. It would be interesting to give an
explicit answer to the question of existence of non-periodic tropical recurrent sequences satisfying $a$.

In section~\ref{four} we design an algorithm which tests for vector
$a=(a_0,\dots,a_n)$ with integer coordinates $a_i\in \ZZ,\, 0\le i\le n$ whether there exists a non-periodic tropical recurrent sequence satisfying $a$.

In section~\ref{five} we introduce a tropical entropy $H(a)$ and a tropical minimal entropy $h(a)$, they fulfil inequalities $0\le h(a) \le H(a)$. We provide an upper bound on $H(a)$ and calculate $h(a)$ and $H(a)$ for some examples of vectors $a$.

In section~\ref{six} we extend the concepts of the tropical (respectively, minimal) entropy to
tropical multivariable recurrent sequences. Again we provide an upper bound on $H(a)$ and calculate $H(a)$ for an example of vector $a$.

Tropical recurrent sequences extend {\it recognizable sequences} \cite{BR}
$z_k=cB^ke,\, 0\le k\in \ZZ$ where $B$ is a square matrix and $c,\, e$ are vectors over a tropical semiring. Observe that recognizable sequences fulfil a tropical linear recurrence due to a tropical version of Cayley-Hamilton theorem \cite{AGG}. Sch\"utzenberger's theorem (see \cite{BR}) characterizes  recognizable sequences by means of rationality of the generating series $\sum_{0\le k < \infty} z_kX^k$. Recognizable sequences are ultimately periodic (so, become periodic for sufficiently large $k$) (see e.~g, \cite{KBR}, \cite{Gaubert}) which differs from tropical recurrent sequences. A related work is \cite{CMQV} where a subclass of non-decreasing recognizable sequences was introduced, and it plays a role in the modelling of timed discrete event systems \cite{BCOQ}.

\section{Tropical recurrent sequences satisfying a vector with Newton graph located on a line
%lying on a single bounded
%edge of Newton polygon
}\label{one}

\begin{definition}\label{minimal}
We say that a tropical recurrent sequence $y=\{y_j\}_{j\in \ZZ}$ satisfying vector $a$ is {\it minimal} if for any $j\in \ZZ$ there exists $k\in \ZZ,\, j-n\le k \le j$ such that $a_{j-k}+y_j=\min_{0\le i\le n} \{a_i+y_{i+k}\}$. In other words,
one can not diminish $y_j$ keeping all the rest $y_l, \, l\neq j$.
%Throughout
%the paper (except for sections~\ref{five}, \ref{six}) we consider only minimal tropical recurrent sequences.
\end{definition}

One can plot $a_i$ as point $(i,a_i)$ on the plane, respectively, $y_j$ as $(j,y_j)$.

\begin{remark}\label{exist}
Making a suitable affine
transformation of the plane one can assume w.l.o.g. that Newton polygon $P(a)$ contains an edge situated on the abscissas axis. Hence $\{y_i=0\, : \, i\in \ZZ\}$ is a tropical recurrent minimal sequence satisfying $a$.
\end{remark}

In this section we study the case when all the points $(i,a_i)$ (so, Newton graph of $a$) are located
on a single bounded edge of Newton polygon $P(a)$. Making a suitable affine
transformation of the plane (cf. Remark~\ref{exist}) one can assume w.l.o.g. that this edge is situated
on the abscissas axis, so either $a_i=0$ or $a_i=\infty$ for all $i\in \ZZ$. In particular, $a_0=a_n=0$.

Note that if $y^{(1)},\, y^{(2)}$ are two tropical recurrent sequences satisfying $a$ then $\min\{b_1+y^{(1)},\, b_2+y^{(2)}\}$ is also the tropical recurrent sequence satisfying $a$, where $b_1,\, b_2 \in \RR$.

There is always the trivial infinite tropical recurrent sequence $\{y_j=\infty,\, j\in \ZZ\}$, so we suppose that  tropical recurrent sequences under consideration differ from the infinite one.

\begin{proposition}\label{one_edge}
Let for vector $a=(a_0,\dots,a_n)$ hold $a_i=0$ for all finite $a_i$. Then all
the tropical recurrent minimal sequences satisfying $a$ are periodic iff all $i$ for
which $a_i=0$ form an arithmetic progression. If $d$ is the difference of this progression then every tropical recurrent minimal sequence satisfying $a$ is periodic with the period $d$.
\end{proposition}

{\bf Proof}. First assume that  set $S:=\{i\in \ZZ: \, a_i=0\}$ does not
form an arithmetic progression. We claim that for any $0<b\in \RR \cup \{\infty\}$ and $k\in \ZZ$ the tropical recurrent sequence $y_i=b$ for $i-k\in S$ and $y_i=0$ otherwise,
satisfies $a$. Suppose the contrary. To simplify notations assume w.l.o.g.
%$$\displaystyle{\phi_k({\bf t})_{k'}=\begin{cases}
%t_{k'} & \text{if $k'<k$}\\
%\sigma_{k,k}({\bf t}) & \text{if $k'=k$}\\
%t_{k'}\cdot \sigma_{k,k'}({\bf t})^{-a_{i_k,i_{k'}}} & \text{if $k'>k$, $i_{k'}\ne i_k$}\\
%\frac{t_{k'}}{\sigma_{k,k'}({\bf t})\cdot
%(t_{k'}+\sigma_{k,k'}({\bf t}))} & \text{if $k'>k$, $i_{k'}=i_k$}\\
%\end{cases}}\ ,$$
%$$\displaystyle{y_i =\begin{cases}
%$b$ & \text{for $i-k\in S$}\\
%$0$ & \text{otherwise}
%\end{cases}}\ $$
%\noindent satisfies $a$
that $k=0$. Then there exists $s \in \ZZ$ such that minimum
$\min_{0\le i\le n} \{a_i+y_{i+s}\}$ is attained once. When $s=0$ we have $a_j+y_j=
b=\min _{0\le i\le n} \{a_i+y_i\}$ for each $j\in S$. In particular, $a_0+y_0=a_n+a_n=b$,
and we get a contradiction with that the minimum is attained once.
When $s>0$ we have
$0=a_n+y_{n+s}=\min_{0\le i\le n}\{a_i+y_{i+s}\}$. Since the minimum is attained once,
we conclude that
%According to the supposition
for each $0\le i<n,\, i\in S$ it holds $a_i+y_{i+s}=b$, therefore $S$ forms
an arithmetic progression with the difference $s$, we get a contradiction.
A similar argument works when $s<0$, in this case $0=a_0+y_s=\min_{0\le i\le n}
\{a_i+y_{i+s}\}$. The claim is proved.

We now show that the sequence $y$ is minimal (see Definition~\ref{minimal}). Indeed, for each $i\in S$
we have $b=a_i+y_i=\min_{0\le l\le n}\{a_l+y_l\}$. On the other hand, for each
$i\notin S$ we have $0=a_0+y_i=\min_{0\le l\le n}\{a_l+y_{i+l}\}$ due to the
proved above claim.

Before proving the other implication in Proposition~\ref{one_edge} we make the following remark.

\begin{remark}\label{modify}
Thus, in case when $S$ does not form an arithmetic progression one can modify
the zero tropical recurrent sequence $\{y_i=0,\, i\in \ZZ\}$ increasing
$y_{i+k}$ for $i\in S$ by $b_k=b>0$. Moreover, one can take arbitrary integers
$\dots,k_{-1},k_0,k_1,\dots$ such that $k_{l+1}-k_l>2n,\, l\in \ZZ$, and for
each $k_l$ modify the zero solution by $b_{k_l}>0$ as described above. Thus,
one obtains an uncountable number of modifications (just by choosing $\{k_l\}_{l\in \ZZ}$ regardless of $b_{k_l}$) of the zero tropical recurrent sequence,
satisfying $a$.
\end{remark}

%Coming back to the proof of Proposition~\ref{one_edge},
Let now $S$ form an
arithmetic progression with a difference $d$ (i.~e. $S=\{id\, : \, 0\le i\le m\}$). Let a tropical recurrent minimal sequence $y$ satisfy $a$.
Fix $0\le i_0 <d$ for the time being. To prove the last statement of Proposition~\ref{one_edge} on periodicity it suffices to
verify that  the subsequence $\{y_{id+i_0}\, : \, i\in \ZZ\}$ of $y$ is constant for $i\in \ZZ$.
%constitutes a tropical recurrent sequence
%satisfying $a$. Therefore, one can consider each of these $d$ subsequences instead of $y$, thus assuming that $d=1$.
%To prove the required last statement in the Proposition on periodicity it suffices to show that $y_i$ is constant for $i\in \ZZ$
%when $d=1$.
Denote $c:=\min_{0\le i\le m}\{y_{id+i_0}\}=\min_{0\le i\le m}\{a_{id}+y_{id+i_0}\}$.
%=\min_{0\le i\le m}\{a_i+y_i\}$.
Then this minimum is attained for two different
$0\le i_1<i_2\le n$.
 %such that $d|(i_2-i_1)$.
Observe that for every $i_1<i<i_2$
%with $d|(i-i_1)$
 it holds $y_{id+i_0}=c$ as well due to the minimality of $y$ (see
Definition~\ref{minimal}). In addition, observe that $c=\min_{-1\le i\le m+1} \{y_{id+i_0}\}$. Indeed, if on the contrary $y_{(m+1)d+i_0}<c$ then minimum $\min_{0\le i\le m}\{a_{id}+
y_{(i+1)d+i_0}\}=y_{(m+1)d+i_0}$ is attained only once, which contradicts that $y$ is a
tropical recurrent sequence satisfying $a$. Similarly, one shows that $y_{-d+i_0}\ge c$.
Repeating this argument recursively one deduces that $c=\min_{-\infty < i <\infty} \{y_{id+i_0}\}$.
Considering $c=a_0+y_{i_2d+i_0}=\min_{0\le i\le m} \{a_{id}+y_{(i+i_2)d+i_0}\}$ one gets that there exists $i_2<i_3\le i_2+m$ such that
$y_{i_3d+i_0}=c$. Then as above one obtains that $y_{id+i_0}=c$ for any $i_2<i<i_3$. Similarly, there exists $i_1-m\le i_4 <i_1$ such that $y_{i_4d+i_0}=c$. Hence $y_{id+i_0}=c$ for any $i_4<i<i_1$. Repeating this argument one concludes that $y_{id+i_0}=c$ for any $i\in \ZZ$.
%Hence $y_{i_1+di}=c$ for any $i\in \ZZ$.

This completes the proof that $y$ is periodic with the period $d$. $\Box$

%One can apply the same argument to any arithmetic progression with the difference $d$ and conclude that for any fixed $i_0$ t%he value of $y_{i_0+di}$ is independent
%from $i$. The periodicity of $y$ is established. $\Box$

\section{Tropical recurrent sequences for a Newton polygon with a single edge and period 2}\label{two}

In the previous section we studied the case when points $(i,a_i),\, 0\le i\le n$ from Newton graph
are located on a line. Note that in general, Newton polygon $P(a)$ has two unbounded edges and
several bounded ones.
In the present section we suppose that $P(a)$ has a single bounded edge. Similar
to the previous section, making a suitable affine transformation of the plane one can
assume w.l.o.g. that this edge lies on the abscissas axis. Again as in the
previous section we consider set $S:=\{0\le i\le n:\, a_i=0\}$. In particular,
$0,n \in S$.

\begin{remark}\label{non-progression}
First, consider the case when $S$ does not form an arithmetic progression. Then
similar to Remark~\ref{modify} one can modify the zero
tropical recurrent sequence $y_j=0,\, j\in \ZZ$ by replacing $y_j=b$ for
$j\in S$, while $b>0$ should be taken less than $\min\{a_i: \, a_i>0\}$. Thus,
again one obtains an uncountable number of periodic tropical recurrent minimal
sequences satisfying $a$.
\end{remark}

So, from now on we assume that $S$ forms an arithmetic
progression with a difference $d$. In the present section we study the case $d=2$ and investigate when
there is a non-periodic tropical recurrent minimal sequence satisfying $a$. In particular, in this case $n$ is even and $2i \in S$ for all $0\le 2i\le n$.

\begin{theorem}\label{two_edges}
Let Newton polygon $P(a)$ have a single bounded edge (assuming w.l.o.g. that this edge is located on the abscissas axis),
and the points $S$ from Newton graph of $a$ on this edge form an arithmetic progression with the
difference 2 starting with the point 0 and finishing with the point $n$. Then any tropical recurrent minimal sequence satisfying $a$ is periodic
iff all $a_i$ with odd $i$ are equal. In the latter case any tropical recurrent minimal
sequence satisfying $a$, has period 2.
\end{theorem}

\begin{remark}\label{parallel}
In other words, under the conditions of the Theorem any tropical recurrent minimal sequence is periodic
iff points $(i,a_i)$ from Newton graph of $a$ are located on two parallel lines: one for even $0\le i\le n$ and the second for odd $1\le i\le n-1$.
\end{remark}

{\bf Proof of the theorem}. First consider the case when not all $a_i$ with odd $i$ are equal. Denote $c:=\min\{a_i:\, \mbox{odd}\, i\}>0,\, e:=\min\{a_i:\, \mbox{odd}\, i,\, a_i>c\}>c$ and $C:=\{i:\, a_i=c\}$. Take a periodic (with the period 2) tropical recurrent sequence $y^{(0)}:=\{y^{(0)}_{2i}=c,\, y^{(0)}_{2i+1}=0,\, i\in \ZZ\}$ satisfying $a$. Let us modify it (denote the modified sequence by $y:=\{y_i\:\, i\in \ZZ\}$) putting

$\bullet$ $y_{2i}:=e,\, 0\le 2i\le n$;

$\bullet$ $y_{2i+1}:=e-c,\, 2i+1\in C, \, 1\le 2i+1 \le n-1$;

$\bullet$ $y_{2i+1}=0,\, 2i+1\notin C,\, 1\le 2i+1 \le n-1$,

\noindent while keeping the rest of
the coordinates unchanged.

Let us verify that $y$ satisfies $a$. For any odd $k\ge 3$ minimum $0=a_n+y_{n+k}=a_{n-2}+y_{n-2+k}=\min_{0\le i\le n} \{a_i+y_{i+k}\}$ is attained at least twice. Similarly, for any odd $k\le -3$ minimum $0=a_0+y_k=a_2+y_{k+2}=\min_{0\le i\le n} \{a_i+a_{i+k}\}$ is  attained at least twice as well. For $k=\pm 1$ minimum $\min_{0\le i\le n} \{a_i+y_{i+k}\}=0$
is attained at least twice since
 %For any odd $k$ minimum
%$\min_{0\le i\le n}\{a_i+y_{i+k}\}=0$ is attained at least twice since
there is an odd $1\le i\le n-1$ such that $i\notin C$. For $k=0$ minimum $\min_{0\le i\le n}
 \{a_i+y_i\}=e$ is attained (for any even $0\le i\le n$, for any $i\in C$ and for any odd $i\notin C$ such that $a_i=e$) at least twice as well. For an even $k\neq 0$ minimum  $\min_{0\le i\le n}\{a_i+y_{i+k}\}=c$.
For an even $k\ge 4$ this minimum is attained at least twice for $c=a_n+y_{n+k}=a_{n-2}+y_{n-2+k}=\min_{0\le i\le n}
\{a_i+y_{i+k}\}$. Similarly, for an even $k\le -4$ this minimum is also attained at least twice for $c=a_0+y_k=a_2+y_{k+2}=\min_{0\le i\le n} \{a_i+y_{i+k}\}$.
%For an even $|k|\ge 4$ this minimum is
%attained at least twice for even coordinates $i$.
Consider the case $k=2$
(the case $k=-2$ can be  considered in a similar way). We have $a_n+y_{n+2}=0+c$.
We claim that there is an odd $1\le 2i-1\le n-1$ such that $a_{2i-1}+y_{2i+1}=c$.
Assume that there exists $2i-1\in C,\, 1\le 2i-1\le n-3$ for which $2i+1\notin C$. Then $a_{2i-1}=c,\, y_{2i+1}=0$, which proves the claim under the assumption.
If such $2i-1\in C$ does not exist then $n-1\in C$, hence $a_{n-1}+y_{n+1}=c+0$,
which proves the claim.

One can check the minimality of $y$ (see Definition~\ref{minimal}) for any (odd)
$2i+1$ such that $y_{2i+1}=0$ with the help of an appropriate even $k$. For an
odd $2i+1\in C$ (in this case $y_{2i+1}=e-c$) one uses $k=0$. Also $k=0$ is used
for even $0\le 2i\le n$ (in this case $y_{2i}=e$). For an even $|2i|\ge n+2$ an
appropriate odd $k$ is involved.

Actually, one can take any $0<q<d-c$ and modify $y^{(0)}$ putting

$\bullet$ $y_{2i}=c+q,\, 0\le 2i\le n$;

$\bullet$ $y_{2i+1}=q,\, 2i+1\in C, \, 1\le 2i+1 \le n-1$;

$\bullet$ $y_{2i+1}=0,\, 2i+1\notin C,\, 1\le 2i+1 \le n-1$,

\noindent while keeping the rest of the coordinates unchanged.

Again as in the proof of Proposition~\ref{one_edge} one can modify $y^{(0)}$ changing
$y_{k_l},y_{k_l+1},\dots,y_{k_l+n}$ as described above for integers
$\dots< k_{-1}<k_0<k_1<\dots$ such that $k_{l+1}-k_l>2n$ for all integers $l$.
Thus, there is an uncountable number of non-periodic tropical recurrent minimal sequences satisfying $a$. \vspace{2mm}

Now we consider $a$ such that $a_{2i}=0,\, 0\le 2i\le n;\, a_{2i+1}=c>0,\, 1\le 2i+1\le n-1$ and prove that any tropical recurrent minimal sequence $y$ satisfying $a$, is
periodic with the period 2.

Denote $b:=\min_{-n\le i\le n}\{y_i\}$. Similar to the proof of Proposition~\ref{one_edge} one can deduce that $b=\min_{-n-1\le i\le n+1}\{y_i\}$ and further by recursion that $b=\min_{-\infty < i < \infty} \{y_i\}$.  Let $y_i=b$ for some
even $-n\le i\le n$ (an odd $i$ can be considered in a similar way). Denote $B:=
\{2i:\, y_{2i}=b\}$. For any pair of adjacent elements $2i_1<2i_2$ of $B$ one has
$2(i_2-i_1)\le n$, because otherwise, minimum $\min_{0\le i\le n}\{a_i+y_{i+2i_1}\}=b$
is attained only once for $i=0$. Therefore, for every even $k$ we have
$\min_{0\le i\le n} \{a_i+y_{i+k}\}=b$, and for every odd $l$ we have $\min_{0\le i\le n} \{a_i+y_{i+l}\}\le b+c$, in addition $y_l\le b+c$ due to the minimality of $y$ (see
Definition~\ref{minimal}). Again due to the minimality $y_k=b$ for every even $k$.

Denote $p:=\min\{y_l:\, -n\le l\le n,\, \mbox{odd}\, l\} \le b+c$. If $p=b+c$ then one deduces that $y_l=b+c$ for any odd $l$ arguing as above by recursion on $|l|$, thus
$y$ is periodic with the period 2. So, assume that $p<b+c$. Arguing as above we
conclude that $p=\min\{y_l:\, \mbox{odd}\, l\in \ZZ\}$ and that for all odd integers
$\dots <l_{-1}<l_0<l_1<\dots$ for which $y_{l_j}=p$, we have $l_{j+1}-l_j\le n$.
Hence due to the minimality of $y$ we get that $y_l=p$ for any odd $l$. Thus, $y$
is periodic with the period 2. Theorem is proved. $\Box$

\section{Newton polygon with a single edge and period greater than 2}\label{three}

So far, we considered vector  $a$ such that its Newton polygon $P(a)$ had
a single bounded edge (recall that w.l.o.g. one can assume that this edge is
situated on the abscissas axis). Moreover, one can suppose that the set $S=\{i\, : \, a_i=0\}$ for the points from Newton graph
on this edge forms an arithmetic progression with a difference $d$ (otherwise, as
we have shown above, there would be a non-periodic tropical recurrent minimal sequence
satisfying $a$). We have given a complete answer to the question of existence of
non-periodic minimal sequences for $d=1,\, 2$. In the present section we provide examples
for $d=3$ which demonstrate that the answer is more complicated in this case and depends not
only on the properties to be arithmetic progressions as for $d=1,\, 2$. It would
be interesting to give a complete answer for $d\ge 3$.

%Thus, we take vector $a: \, a_0=a_3=0,\, a_1=b,\, a_2=c;\, b,c>0, b\neq c$.

\begin{proposition}\label{degree_3}
Let vector $a: \, a_0=a_3=0,\, a_1=b,\, a_2=c;\, b,c>0, b\neq c$.
Any tropical recurrent minimal sequence satisfying $a$ is periodic iff either $b<c\le 2b$
or $c<b\le 2c$. In this case any minimal sequence is periodic with the period 3.
\end{proposition}

{\bf Proof}. First consider the case $c>2b$ (the case $b>2c$ can be studied in a
similar way). Take the following periodic tropical recurrent minimal sequence satisfying $a$:
$y^{(0)}_{3i}:=0;\, y^{(0)}_{3i+1}=2b;\, y^{(0)}_{3i+2}:=b;\, i\in \ZZ$. Consider a real $0<e\le c-2b$.
We modify $y^{(0)}$ resulting in a non-periodic tropical recurrent minimal sequence $y$ satisfying $a$:
$y_1=2b+e,\, y_2=b+e,\, y_4=2b+e$ and keeping the rest of the coordinates of $y^{(0)}$ unchanged.

Similar to the proofs of Proposition~\ref{one_edge} and Theorem~\ref{two_edges} one can choose integers
$\dots < k_{-1}<k_0<k_1<\dots$ and reals $0<e_l\le c-2b$ such that $k_{l+1}-k_l\ge 3$, and modify $y^{(0)}$
putting $y_{3k_l+1}=2b+e_l;\, y_{3k_l+2}=b+e_l;\, y_{3k_l+4}=2b+e_l$ for all integers $l$. Thus, one achieves
a non-countable number of non-periodic tropical recurrent minimal sequences satisfying $a$.  \vspace{2mm}

Now we study $a$ fulfilling inequalities $b<c\le 2b$ (the case of inequalities
$c<b\le 2c$ is considered in a similar way). Let a tropical recurrent minimal sequence $y$  satisfy $a$. Denote $q:=\min_{-3\le i\le 3}\{y_i\}$. Arguing as in the proofs of Proposition~\ref{one_edge} and Theorem~\ref{two_edges}, we conclude that $q=\min_{-4\le i\le 4}\{y_i\}$, and continuing this argument we get by induction that
$q=\min_{-\infty<i<\infty}\{y_i\}$. One can assume w.l.o.g. that $y_{3i_0}=q$ for
some integer $i_0$. Since minimum $\min_{0\le i \le 3} \{a_i+y_{3i_0+i}\}$ is
attained at least twice, we deduce that $y_{3i_0+3}=q$. Continuing in this way,
we deduce that $y_{3j}=q$ for every integer $j$.

Now we consider coordinates for an arithmetic progression $\{y_{3i+2}:\, i\in \ZZ\}$. Denote $r:=\min_{0\le i\le 1} \{y_{3i+2}\}$. Hence $r\le q+b$ since minimum
$\min_{0\le i\le 3}\{a_i+y_{i+2}\} \le a_1+y_3=b+q<a_2+q\le c+y_4$ should be attained at least twice. Arguing as above, we deduce that $r=\min_{-\infty < i < \infty}\{y_{3i+2}\}$. If $r<q+b$ then since $\min_{0\le i\le 3} \{a_i+y_{3j+i+2}\} \le a_0+r=r<
b+q=a_1+y_{3j+3}<c+q\le a_2+y_{3j+4}$ for any integer $j$, we conclude that
$y_{3i+2}=r$ for every integer $i$. Now let $r=q+b$. Since $\min_{0\le i\le 3}
\{a_i+y_{3j+i+1}\}\le a_2+q=c+q\le 2b+q=b+r\le a_1+y_{3j+2}$, we get that if
$y_{3j+2}>r$ then $a_1+y_{3j+2}$ does not attain minimum in $\min_{0\le i\le 3}
\{a_i+y_{3j+i+1}\}$, hence $a_0+y_{3j+2}=y_{3j+2}$ attains the minimum in
$\min_{0\le i\le 3}\{a_i+y_{3j+i+2}\}\le a_1+y_{3j+3}=b+q=r$ according to the
minimality of $y$ at $y_{3j+2}$ (cf. Definition~\ref{minimal}). Thus, $y_{3i+2}=r$
for every integer $i$.

Finally, we consider coordinates for an arithmetic progression $\{y_{3i+1}: \,
i\in \ZZ\}$. Due to the minimality of $y$ one deduces that $y_{3i+1}\le t_0:=
\min\{b+r,\, c+q\}$ for every integer $i$. Denote $t:=\min_{0\le i\le 1}\{y_{3i+1}\}\le t_0$. Arguing as above, we conclude that $t=\min_{-\infty < i < \infty}\{y_{3i+1}\}$. When $t=t_0$, we have obviously, $y_{3i+1}=t$ for every integer $i$. When
$t<t_0$, arguing as above we also show that $y_{3i+2}=t$ for every integer $i$.
This completes the proof of the  Proposition. $\Box$

\begin{remark}\label{multi_edge}
When $P(a)$ has several edges $g_1,\dots,g_k$ with slopes $\sigma_1<\cdots <\sigma_k$, respectively,
then as a tropical recurrent sequence satisfying $a$ one can take points
$\{(i,\, y_i): \, i\in \ZZ\}$ lying on the edges of an (infinite in both directions) convex polygon having edges $g_k',\, g_{k-1}',\dots,g_1'$ with the slopes
$-\sigma_k,-\sigma_{k-1},\dots,-\sigma_1$, respectively, such that  $g_j'$ is not shorter than $g_j,\, 1\le j\le k$.

Conversely, in section 4 \cite{G} it is proved, in fact, that if for each $0\le i\le n$ point $(i,\, a_i)$ lies on the boundary of $P(a)$ then any tropical recurrent sequence satisfying $a$ has the described form.
\end{remark}

\section{Algorithm testing existence of a non-periodic tropical recurrent sequence} \label{four}

Let $a=(a_0,\dots,a_n)\in \ZZ$ be a vector with integer coordinates whose Newton polygon has a single bounded edge which is located on the abscissas axis.

\begin{definition}\label{cut}
We say that a {\it vector $z=(z_0,\dots,z_N)\in \RR^{N+1}$ satisfies $a$ and is minimal} if
for each $0\le k\le N-n$ minimum $\min_{0\le i\le n} \{a_i+z_{i+k}\}$ is
attained at least twice and for each $n\le j\le N-n$ there exists $j-n\le k\le j$ such that
$a_{j-k}+z_j=\min_{0\le i\le n} \{a_i+z_{k+i}\}$
(cf. the definition of a tropical recurrent sequence and Definition~\ref{minimal}).
\end{definition}

In this section we prove the following theorem.

\begin{theorem}\label{algorithm}
There is an algorithm which for a vector $a=(a_0,\dots,a_n),\, 0\le a_i\le M,\, 0\le i\le n,\, a_0=a_n=0$ tests whether there exists a non-periodic tropical recurrent minimal sequence satisfying $a$. The complexity of the algorithm does not exceed $(Mn)^{O(n)}$.
\end{theorem}

{\bf Proof}. 
%It suffices to 
Consider a tropical recurrent minimal sequence $y=\{y_i\in \RR,\, i\in \ZZ\}$ satisfying $a$.
%with integer coordinates $y_i$.
Denote
$q:=\min_{-n\le i\le n} \{y_i\}$. Replacing $y_i$ by $y_i-q,\, i\in \ZZ$ one can assume w.l.o.g. that $q=0$. Arguing as in the proofs of Proposition~\ref{one_edge} and of Theorem~\ref{two_edges} above, we conclude that
$\min_{-\infty < i < \infty} \{y_i\} =0$.

Denote $S:=\{i\in \ZZ:\, y_i=0\}$. Again arguing as in the proofs of Proposition~\ref{one_edge} and of Theorem~\ref{two_edges} above, we deduce that for any pair of adjacent element $i_1<i_2$ of $S$ inequality $i_2-i_1\le n$ holds.

For each $j\in \ZZ$ due to the minimality of $y$ (see Definition~\ref{minimal}) there exists an integer
$k,\, j-n\le k \le j$ such that $a_{j-k}+y_j=\min_{0\le i\le n}\{a_i+y_{k+i}\}$. On the other hand, there is
$0\le l\le n$ such that $y_{k+l}=0$, hence $y_j\le a_l$. Thus, $y_j\le M,\, j\in \ZZ$.

%A  vector $z=(z_0,\dots,z_N)\in \{0,\dots,M\}^{N+1}$ we treat
%satisfies $a$ if
%for each $0\le k\le N-n$ minimum $\min_{0\le i\le n} \{a_i+z_{i+k}\}$ is
%attained at least twice and for each $n\le j\le N-n$ there exists $j-n\le k\le j$ such that
%$a_{j-k}+z_j=\min_{0\le i\le n} \{a_i+z_{k+i}\}$
%(cf. the definition of a tropical recurrent sequence and Remark~\ref{minimal}).
%Treating $z$
Denote $V:=\{j+i/(2n+2)\, : \, 0\le j<M,\, 0\le i\le 2n+2\}$. We consider words $z=(z_0,\dots,z_N)\in V^{N+1}$ 
in the alphabet $V$ agreeing that a subword of a word consists of its consecutive letters. We correspond to a sequence $y=\dots y_{-1}y_0y_1\dots$ the sequence $Y:=Y(y)$ in the alphabet $W:=V^{2n+1}$ consisting of all the subwords of $y$ of the length $2n+1$, so that $(y_{i+1}\dots y_{i+2n+1})\in W$ is the next letter in $Y$ after $(y_i\dots y_{i+2n})\in W$ for all $i\in \ZZ$. Later on in the proof of Theorem~\ref{algorithm} we consider arbitrary sequences over the reals (rather than over $V$).

Construct a directed graph $G_0:=G_0(a)$ whose vertices are elements of $W$ which satisfy $a$ and being minimal (see Definition~\ref{cut}). There is an arrow in $G_0$ from a vertex $(z_0\dots z_{2n})\in W$ to a vertex $(z_0'\dots z_{2n}')\in W$ iff $z_{i+1}=z_i',\, 0\le i\le 2n-1$ (in particular, loops are admitted). In other words, vector $z_0\dots z_{2n}z_{2n}'=z_0z_0'\dots z_{2n}' \in V^{2n+2}$ satisfies $a$ and is minimal.

\begin{lemma}\label{graph}
(i) Let a tropical recurrent minimal sequence $y$ over $V$ satisfy $a$. Then for any letter from $W$ of the sequence $Y=Y(y)$ there is an arrow in $G_0$ from this letter to the next letter in $Y$. Thus, $Y$ provides an infinite in both directions path $p(Y)$ in $G_0$;

(ii) conversely, for any infinite in both directions path $p$ in $G_0$ there is a unique tropical recurrent minimal sequence $y$ over $V$ satisfying $a$ such that $p=p(Y(y))$;

(iii) sequence $y$ is periodic with a period $d$ iff $Y$ is periodic with the period $d$ as well.
\end{lemma}

{\bf Proof}. The only non-trivial statement is (ii). To yield $y=\dots y_{-1}y_0y_1\dots$ take consecutively the first components (being elements of $V$) of the letters (from $W$) of sequence $Y$ produced from $p$. Then $Y=Y(y)$ and $p=p(Y(y))$. To prove minimality of $y$ take its arbitrary letter $y_j$ and make use of that the word $Y_j:=(y_{j-n}\dots y_j\dots y_{j+n}) \in W$ is a minimal vector satisfying $a$ since $Y_j$ is a vertex of $G_0$, and thus $a_{j-k}+y_j=\min_{0\le i\le n} \{a_i+y_{k+i}\}$ for some $j-n\le k\le j$ (see Definition~\ref{cut}). The proof of that $y$ satisfies $a$ is similar (and even easier). $\Box$ \vspace{1mm}

Observe that vertices of $G_0$ without in-coming or out-coming arrows do not appear in infinite in both directions paths, that is why we remove such vertices successively starting from $G_0$, while it is possible to remove, and the resulting (directed) graph denote by $G:=G(a)$.
Remark~\ref{exist} implies non-emptyness of $G$.

%Theorem~\ref{algorithm} follows from the next lemma.

\begin{lemma}\label{cycle}
All the tropical recurrent minimal sequences over $V$ satisfying $a$ are periodic iff $G(a)$ is a disjoint union of simple cycles.
\end{lemma}

{\bf Proof}. It suffices to show that if one of the connected components $K$ of $G=G(a)$ (to define connected components we consider $G$ as an undirected graph) is not a simple cycle then there exists a non-periodic path in $G$ (see Lemma~\ref{graph}). To show this, treat $G$ as a partial order on its vertices. Then the vertices of $G$ are partitioned into classes of equivalence (two vertices belong to the same class if there exists a cycle in $G$ containing both vertices). The classes of equivalence are partially ordered (this order is induced by $G$).

Since $G$ contains no vertices without in-coming or out-coming arrows, every minimal or maximal (with respect to the mentioned above partial order) class of equivalence contains a cycle with at least one arrow. Therefore,
either there is (in $K$) a class of equivalence being not a simple cycle or there are at least two classes of equivalence. First, let there be a class of equivalence $C$ in $K$ being not a simple cycle. Then $C$ contains suitable three vertices $v_0,\, v_1,\, v_2$ and two arrows $(v_0,\, v_1),\, (v_0,v_2)$. There exists a non-periodic path passing through the vertices of $C$, and each time passing $v_0$ the path chooses either $v_1$ or $v_2$ as the next vertex in an arbitrary non-periodic manner.

Second, let there be at least two classes of equivalence in $K$. Then there exists a path in $G$ of the form $v_0v_1\dots v_k(=v_0)u_1\dots u_lv_0'v_1'\dots v_q'(=v_0')$ such that the cycle $v_0v_1\dots v_k(=v_0)$ belongs to one class of equivalence, while the cycle $v_0'v_1'\dots v_q'(=v_0')$ belongs to another class of equivalence. As a non-periodic path in $G$ one can take a repetition of an infinite number of times of cycle $v_0v_1\dots v_k(=v_0)$ (for negative moments of time till 0), subsequently $u_1\dots u_l$ (for the moments $1,\dots l$), and
subsequently a repetition of an infinite number times of cycle $v_0'v_1'\dots v_q'(=v_0')$ (starting with the moment $l+1$). $\Box$ \vspace{1mm}

Lemmas~\ref{graph}, \ref{cycle} provide an algorithm to test whether there exists a non-periodic tropical recurrent minimal sequence over $V$ satisfying $a$ (cf. Theorem~\ref{algorithm}).

\begin{remark}
Let us call a sequence $y=\{y_i\, : \, i\in \ZZ\}$ {\it stable periodic} if for suitable $d$ we have $y_{i+d}=y_i$ for all sufficiently big absolute values $|i|$. One can deduce similar to the proof of Lemma~\ref{cycle} that all the tropical recurrent minimal sequences over $V$ satisfying $a$ are stable periodic iff all the classes of equivalence of vertices of $G(a)$ (defined in the proof of Lemma~\ref{cycle}) are simple cycles.
\end{remark}

Denote by $\langle e \rangle:= e-\lfloor e \rfloor$ the fractional part of $e\in \RR$. We will complete the proof of Theorem~\ref{algorithm} by reducing it to Lemmas~\ref{graph}, \ref{cycle} relying on the following lemma.

\begin{lemma}\label{equalizing}
Let $0=c_0<c_1<\cdots <c_m<c_{m+1}=1,\, 0=e_0<e_1<\cdots <e_m<e_{m+1}=1$ be a pair of real increasing sequences. Assume that a tropical recurrent minimal sequence $y=\{y_i\, : \, i\in \ZZ\}$ satisfies $a$. For each $i\in \ZZ$ take a unique $0\le j\le m$ such that $c_j\le \langle y_i\rangle <c_{j+1}$ and replace $y_i$ by $x_i:=\lfloor y_i\rfloor +e_j$. Then $x:=\{x_i\, : \, i\in \ZZ\}$ is also a tropical recurrent minimal sequence satisfying $a$.
\end{lemma}

{\bf Proof}. Lemma follows from the claim that the described replacement keeps all the non-strict inequalities appearing in Definition~\ref{minimal}. Namely, for $0\le i_1<i_2\le n$ and $k\in \ZZ$ let $y_k+a_{i_1}\le y_{k+i_2-i_1}+a_{i_2}$ (the opposite inequality $y_k+a_{i_1}\ge y_{k+i_2-i_1}+a_{i_2}$ can be considered in a similar manner), then we claim that $x_k+a_{i_1}\le x_{k+i_2-i_1}+a_{i_2}$.

Indeed, first suppose that there exists an integer $t\in \ZZ$ such that $y_k+a_{i_1}\le t\le y_{k+i_2-i_1}+a_{i_2}$. Then the similar inequalities $x_k+a_{i_1}\le t\le x_{k+i_2-i_1}+a_{i_2}$ hold.

Otherwise, $t<y_k+a_{i_1}\le y_{k+i_2-i_1}+a_{i_2}<t+1$ for a suitable integer $t\in \ZZ$. Hence $t=\lfloor y_k\rfloor +a_{i_1}=\lfloor y_{k+i_2-i_1}\rfloor +a_{i_2}$. There exist unique $0\le j_1,j_2\le m$ such that $c_{j_1}\le \langle y_k\rangle < c_{j_1+1},\, c_{j_2}\le \langle y_{k+i_2-i_1}\rangle <c_{j_2+1}$. Then $j_1\le j_2$. Therefore, 
$$x_k+a_{i_1}=\lfloor y_k \rfloor +e_{j_1}+a_{i_1}=t+e_{j_1}\le t+e_{j_2}= \lfloor y_{k+i_2-i_1} \rfloor +e_{j_2}+a_{i_2}=x_{k+i_2-i_1}+a_{i_2}.$$
\noindent The claim is proved. $\Box$ 
%\vspace{1mm}

\begin{lemma}
If there exists a non-periodic tropical recurrent minimal sequence over the reals then there exists a non-periodic tropical recurrent sequence over $V$.
\end{lemma}

{\bf Proof}. Let a non-periodic tropical recurrent minimal sequence $y=\{y_i\in \RR\, :\, i\in \ZZ\}$ satisfy $a$. First assume that among all the fractional parts $0=c_0<c_1< \cdots < c_m<c_{m+1}=1$ of $\langle y_i \rangle,\, i\in \ZZ$ there are at most $m\le 2n+1$ different. Then one can apply Lemma~\ref{equalizing} and obtain a tropical recurrent minimal sequence $x:=\{x_i\, :\, i\in \ZZ\}$ with the fractional parts $0=e_0<e_1< \cdots < e_m<e_{m+1}=1$ where $e_j=j/(2n+2),\, 0\le j\le m$. Observe that $x$ is also non-periodic since for any pair of integers $i_1<i_2$ we have $y_{i_1}=y_{i_2}$ iff $x_{i_1}=x_{i_2}$, taking into account that if $y_i=\lfloor y_i\rfloor +c_j$ for some $0\le j\le m$ then $x_i=\lfloor y_i\rfloor +e_j$. Note that $x$ is a sequence in the alphabet $V$ because $0\le y_i\le M,\, i\in \ZZ$ as it was shown at the beginning of the proof of Theorem~\ref{algorithm}.

Now we assume that there are more than $2n+1$ different fractional parts among $\langle y_i \rangle,\, i\in \ZZ$. Pick $i_0\in \ZZ$ fulfilling the property that to the right of the subword $y_{i_0}\cdots y_{i_0+2n}$ there exists $i>i_0+2n$ such that the fractional part $\langle y_i \rangle$ is not among fractional parts $\langle y_{i_0} \rangle,\dots,\langle y_{i_0+2n} \rangle$. Take the minimal such $i$. Denote by $0=c_0<c_1< \cdots < c_m<c_{m+1}=1$ the (ordered) fractional parts $\langle y_s \rangle,\, i-2n-1\le s\le i$, hence $m\le 2n+1$. Let $\langle y_i \rangle = c_{j_0}$ for some $1\le j_0\le m$ (note that $\langle y_i \rangle \neq 0$ since among $\langle y_s \rangle,\, i-2n-1\le s<i$ at least twice occurs $0$ because for adjacent $s_1<s_2$ with $y_{s_1}=y_{s_2}=0$ we have $s_2-s_1\le n$ as it was shown above at the beginning of the proof of Theorem~\ref{algorithm}).

We apply Lemma~\ref{equalizing} to $y$ twice: first time with $0=c_0<c_1< \cdots < c_m<c_{m+1}=1$ and $e_j=j/(2n+2),\, 0\le j\le m$, the resulting sequence denote by $x:=\{x_s\, :\, s\in \ZZ\}$. Second time we apply Lemma~\ref{equalizing} to $y$ with $0=c_0<c_1< \cdots c_{j_0-1} < c_{j_0+1}<\cdots <c_m<c_{m+1}=1$ and $e_j=j/(2n+2),\, 0\le j\le m,\, j\neq j_0$, the resulting sequence denote by $x':=\{x'_s\, :\, s\in \ZZ\}$. Both sequences $x,x'$ are over the alphabet $V$. Therefore, the vertex $w:=x_{i-2n-1}x_{i-2n}\dots x_{i-1}=x'_{i-2n-1}x'_{i-2n}\dots x'_{i-1}\in W$ belongs to the graph $G$ (cf. Lemma~\ref{cycle}), and there are two arrows in $G$ with the same left end-point $w$ and with the right end-points $x_{i-2n}\dots x_{i-1}x_i \neq x'_{i-2n}\dots x'_{i-1}x'_i$, respectively, being different since $\langle x_i \rangle = e_{j_0}> \langle x'_i \rangle = e_{j_0-1}$. Hence, Lemma~\ref{cycle} implies that there exists a non-periodic tropical recurrent minimal sequence over the alphabet $V$ satisfying $a$. Lemma is proved. $\Box$ \vspace{1mm}

Thus, the algorithm required in Theorem~\ref{algorithm} constructs the graph $G$ and tests whether it is a disjoint union of simple cycles that completes the proof of the Theorem. $\Box$ \vspace{1mm}

\begin{remark}
(i) More generally, one can design an algorithm similar to Theorem~\ref{algorithm} for vectors $a=(a_0,\dots,a_n)$ with rational coordinates $a_i\in \QQ,\, 0\le i\le n$ looking for tropical recurrent minimal sequences of the form $y=\{y_i/q:\, i\in \ZZ,\, y_i\in \ZZ\}$ where $q$ being the common denominator of $a_0,\dots,a_n$;

(ii) it would be interesting to design an algorithm similar to Theorem~\ref{algorithm} for vectors $a$ allowing real algebraic and infinite coordinates.
\end{remark}

\section{Tropical entropy}\label{five}

Let $a=(a_0,\dots,a_n)\in \RR^{n+1}$. For $0\le s\in \ZZ$ denote by $D_s\subset \RR^s$ (respectively, $M_s\subset \RR^s$) the set of vectors satisfying $a$ (respectively, satisfying in addition the minimality condition, see Definition~\ref{cut} and section~\ref{four}). Both $D_s$ and $M_s$ are polyhedral complexes \cite{MS}. Denote $d_s:=\dim D_s,\, m_s:=\dim M_s$.

When $i+j=s$ denote by $p:\RR^s \to \RR^i,\, q:\RR^s\to \RR^j$ the projection onto the first $i$ coordinates and respectively, onto the last $j$ coordinates. Since $p(D_s)\subset D_i,\, q(D_s)\subset D_j,\, p(M_s)\subset M_i,\, q(M_s)\subset M_j$ we have $d_{i+j}\le d_i+d_j,\, m_{i+j}\le m_i+m_j$. Therefore, due to Fekete's subadditive lemma \cite{S} there exist limits (which coincide as well with the infimums of these sequences)
$$H=H(a):=\lim_{s\to \infty} d_s/s, \quad h=h(a):=\lim_{s\to \infty} m_s/s$$
\noindent which we call the {\it tropical entropy of $a$} (respectively, the
{\it tropical minimal entropy of $a$}). Clearly, $0\le h\le H$.

\begin{proposition}\label{upper}
$H\le 1-1/n$
\end{proposition}

{\bf Proof}. The polyhedral complex $D_s$ is a union of polyhedra such that each of these polyhedra $Q$ satisfies the following conditions. For every $0\le j\le s-n$ there exists a pair $0\le i_1<i_2\le n$ such that $z_{j+i_1}+a_{i_1}=z_{j+i_2}+a_{i_2}=\min_{0\le i\le n} \{z_{j+i}+a_i\}$ for any $(z_1,\dots,z_s)\in D_s$. For $j=i_1+1$ there exists a pair $0\le i_3<i_4\le n$ fulfilling the similar conditions, hence $(i_1+1+i_3)-(j+i_1)\le n$.
Therefore, there are at least $\lfloor s/n\rfloor$ such pairs. Each such pair $(j+i_1,\, j+i_2)$ imposes a linear restriction $z_{j+i_1}-z_{j+i_2}=a_{i_2}-a_{i_1}$ on $Q$.
Thus, $d_s\le s-\lfloor s/n\rfloor$. $\Box$

\begin{example}\label{segment}
Let now vector $a$ be with $a_0=\cdots=a_n=0$. Consider a polyhedron consisting of vectors $z=(z_1,\dots,z_s)\in \RR^s$ such that $z_{(n+1)i+1}=z_{(n+1)i+2}=0$ for every $i$ and the rest of the coordinates being arbitrary real non-negative. Then $z\in D_s$, hence $d_s\ge s(1-2/(n+1))$, thus, $H(a)\ge 1-2/(n+1)$.

Now we prove an upper bound on $H$. Let $y=(y_1,\dots,y_s)\in D_s$. Suppose w.l.o.g. that
$\min_{1\le i\le s} \{y_i\}=0$. Denote $I:=\{1\le i\le s \, : \, y_i=0\}$ and $i_1<i_2<\cdots$ being the consecutive elements of $I$. Observe that $i_{j+2}-i_j\le n$ for every $j$ since otherwise, $\min_{0\le l\le n} \{a_l+y_{l+i_j+1}\}$ is attained only once for $l=i_{j+1}-i_j-1$. Therefore, $y$ belongs to a linear space $\{y_i=0,\, i\in I\}$ with the dimension at most $\lceil s(1-2/(n+1))\rceil$. Thus, $H(a)=1-2/(n+1)$.

Due to Proposition~\ref{one_edge} any vector satisfying $a$ and the condition of the minimality (see Definition~\ref{minimal}) has all equal coordinates, so $m_s=1$, hence $h(a)=0$.
\end{example}

\begin{remark}
There is a gap between an upper bound on $H$ from Proposition~\ref{upper} and the latter example. The conjecture is that $1-2/(n+1)$ is an upper (sharp) bound on $H$.
\end{remark}

We call a vector $a$ {\it regular} if the set $J:=\{i\, |\, a_i<\infty\}$ is an arithmetic progression and each point $(i,\, a_i),\, i\in J$ is a vertex of the Newton polygon $P(a)$.

\begin{theorem}\label{lower}
If a vector $a$ is not regular then $H(a)\ge 1/6$.
\end{theorem}

%\begin{proposition}\label{lower}
%(i) Let a vector $a$ be such that its Newton polygon contains at least three
%points of $a$ on one its edges. Then $H(a)\ge 1/4$;

%(ii) Let Newton polygon $P(a)$ of a vector $a=(a_0,\dots,a_n)$ have a single bounded edge, and this edge  contain only its two endpoints
%among the points of $a$ and $a$ have at least three points. Then $H(a)\ge 1/6$.
%. In addition, $\min_{1\le i\le n-1} \{a_i\} < \infty$. Then $H(a)\ge 1/(2n)$.
%\end{proposition}

{\bf Proof}. First consider the case when at least three points of $a$ lie on a (bounded) edge of $P(a)$.
 Similar to the beginnings of sections~\ref{one},~\ref{two} making
suitable affine transformations one can suppose w.l.o.g. that an edge containing
at least three points of $a$ lies on the abscissas axis. Consider the points of
$a$ located on this edge: $\{(i,\, 0)\, : \, i\in I\}$ where $|I|\ge 3$. One can assume w.l.o.g. that the greatest common divisor $GCD(I)$ of the differences
$i_1-i_2$ of all the pairs of the elements $i_1,\, i_2\in I$ of $I$ equals 1.
Otherwise, one can consider separately all $GCD(I)$ arithmetic progressions with
the difference $GCD(I)$.

Pick any three elements of $I$ not all with the same parity, say $0,\, 2i,\, j$
w.l.o.g. where $i\ge 1$ and $j$ being odd. Consider the following tropical
recurrent sequence $z$ satisfying $a$. For odd indices $2l+1$ we put
$z_{2l+1}=0$. For even indices $2(2qi+t),\, q\in \ZZ,\, 0\le t<i$ we put
$z_{2(2qi+t)}=0$ and $z_{2((2q+1)i+t)}$ we put arbitrarily non-negative. Taking
finite fragments $(z_0,\dots,z_N)$ with growing $N$ we conclude that $H(a)\ge 1/4$ in this case.

Now assume that no edge of $P(a)$ contains a point of $a$ other than two vertices of this edge. Take an edge of $P(a)$ with the vertices $(i,\, a_i),\, (j,\, a_j)$ with the maximal difference $j-i>0$. Again one can suppose w.l.o.g. that these vertices are $(0,\,0),$ and $(n,\, 0)$. There exists $i\in J$ such that $n$ does not divide $i$ since $a$ is not regular. Among such $i$ pick $i_0$ for which $c:=a_{i_0}$ is minimal, then $c>0$. Denote $k:=GCD(n,i_0)$.

%\vspace{2mm}
%(ii) Making a suitable affine transformation of the plane as above, one can ass%ume w.l.o.g. that the single bounded edge of $P(a)$ lies on the abscissas axis. Thus, $a_0=a_n=0$ and $a_i>0,\, 1\le i\le n-1$.

%Denote $c:=\min_{1\le i\le n-1} \{a_i\}>0$ and let $a_{i_0}=c$ for an appropriate $1\le i_0\le n-1$.  Denote $k:=GCD(n,\, i_0)$.
 Consider a sequence $z:=\{z_i\}_{i\in \ZZ}$ such that

$\bullet$ $z_{qn-2ji_0+i}=0$ when $0\le 2j <n/k$;

$\bullet$ $z_{2qn-(2j+1)i_0+i}=c$ when $0< 2j+1<n/k$;

$\bullet$ $z_{(2q+1)n-(2j+1)i_0+i}\ge c$ when $0< 2j+1<n/k$

\noindent for $q\in \ZZ,\,  0\le i< k$. Then $z$ satisfies $a$. Taking
finite fragments $(z_0,\dots,z_N)$ with growing $N$ we conclude that $H(a)\ge 1/4$ for even $n/k$ and $H(a)\ge  \frac{k\lfloor n/(2k) \rfloor}{2n} \ge 1/6$ (the latter inequality becomes an equality when $n/k=3$). $\Box$

\begin{example}
Let vector $a:=(0,\, 1,\, 0)$. Then for any sequence $(z_0,\dots,z_N)$ satisfying $a$ denote $c:=\min\{z_0,\dots,z_N\}$. Let $z_i=c$ for some $0\le i\le N$. For definiteness assume that $i$ is even (the case of an odd $i$ is considered in a similar way). Then $z_l=c$ for any even $0\le l\le N$. If $m:=\min_{0<2j+1\le N} \{z_{2j+1}\}<c+1$ then $z_{2j+1}=m$ for any $0<2j+1\le N$. If $m\ge c+1$ then $m=c+1$. For any odd $0<2j-1\le N-2$ we have that either $z_{2j-1}=c+1$ or $z_{2j+1}=c+1$. Therefore, the number of odd $0<2j+1\le N$ for which $z_{2j+1}>c+1$ does not exceed $\lceil N/4\rceil$, thereby $H(a)\le 1/4$.

On the other hand, from the proof of Theorem~\ref{lower} we conclude that $H(a)\ge 1/4$ (in the notations of the proof of Theorem~\ref{lower} $n=2,\, i_0=1$).
\end{example}

%$\bullet$ $y_i=0$ when $i=ns+l,\, 1\le l\le n-1$;

%$\bullet$ $y_i=c$ when $i=2ns$;

%$\bullet$ $y_i\ge c$ when $i=2ns+n$

%\noindent for $s\in \ZZ$.

%For $k=ns+l, \, 1\le l\le n-1$ minimum $0=a_0+y_k=a_n+y_{n+k}=\min_{0\le i\le n} \{a_i+y_{i+k}\}$ is attained at least t%wice. For $k=2ns$ minimum $c=a_0+y_k=a_{i_0}+y_{i_0+k}=\min_{0\le i\le n} \{a_i+y_{i+k}\}$ is also attained at least t%wice. Finally, for $k=2ns+n$ minimum $c=a_n+y_{k+n}=a_{i_0}+y_{i_0+k}=\min_{0\le i\le n} \{a_i+y_{i+k}\}$ is attained t%wice as well.

%Therefore, $H(a)\ge 1/(2n)$.
%$\Box$

\begin{remark}
If all the points $(i,a_i),\, a_i<\infty,\, 0\le i\le n$ are the vertices of the Newton
polygon $P(a)$ (thus, $a$ being regular) then section 4 \cite{G} implies that $H(a)=0$ (cf. also Remark~\ref{multi_edge}). Our conjecture is that a stronger than in Theorem~\ref{lower} bound $H(a)\ge 1/4$ holds for not regular vector $a$.
% Observe that Proposition~\ref{lower} entails that if $P(a)$ has  a single bounded edge and $a$ contains at least three points then $H(a)\ge 1/6$. How to describe vectors $a$ for which $H(a)=0$?
\end{remark}
%Is the converse statement true? It is true when $P(a)$ has a single bounded edge due to Proposition~\ref{lower};

\begin{corollary}
A vector $a$ is regular iff $H(a)=0$. For non-regular $a$ the inequality $H(a)\ge 1/6$ holds.
\end{corollary}

\begin{example}
Now we give an example of a vector $a=(a_0,\, a_1,\, a_2,\, a_3)$ with a positive $h(a)$. Put $a_0=a_3=0,\, a_2>2a_1>0$ (see Proposition~\ref{degree_3}). Then at the beginning of the proof of Proposition~\ref{degree_3} a family of vectors satisfying $a$ is constructed for an arbitrary $k_1<k_2<\cdots$ with $k_{i+1}-k_i=3$ for all $i$. which demonstrates that $h(a)\ge 1/9$.

On the contrary, if $2a_1\ge a_2 >a_1 >0$ then $h(a)=0$ (cf. the proof of Proposition~\ref{degree_3}).
\end{example}

\begin{remark}
(i) For $a$ with a single bounded edge of its Newton polygon (cf. sections~\ref{one}, \ref{two}, \ref{three}) is it true that $h(a)=0$ iff all the tropical minimal recurrent sequences satisfying $a$ are periodic?

% in case when all the tropical recurrent sequences satisfying $a$ are periodic (cf. Proposition~\ref{one_edge}, Theorem~\ref{two_edges},
 %Proposition~\ref{degree_3}) we have $h(a)=0$. What one can say about a converse statement?

(ii) For a vector $a$ from Example~\ref{segment} it was, actually, established
that when $(n+1)|(s-r),\, 0\le r\le n$, function $d_s=(s-r)(n-1)/(n+1)+r$ in case if $0\le r<n$ and $d_s=(s-n)(n-1)/(n+1)+n-1$ if $r=n$. Is it true that for an arbitrary vector $a$ function $d_s$ is linear with the leading coefficient $H(a)$ for $s$ from each fixed arithmetic progression with the difference $n+1$?
%are $d_s,\, m_s$ linear functions in $s$ (for sufficiently big $s$)?
\end{remark}

\section{Tropical multivariable recurrent sequences}\label{six}

For a vector $a=\{a_I\in \RR\cup \{\infty\}\, :\, I\in \ZZ^m\}$ with a finite number
of $I$ such that $a_I\in \RR$ (the set of such $I$ we call the {\it support} of $a$)
we say that a {\it tropical multivariable recurrent sequence} $\{z_I\in \RR\, :\,
I\in \ZZ^m\}$ {\it satisfies} $a$ if for any $J\in \ZZ^m$ the minimum $\min_I\{a_I+
z_{I+J}\}$ is attained at least for two different $I$. Similar to Definition~\ref{minimal} one can also
define tropical multivariable recurrent {\it minimal}  sequences.

For a parallelepiped in the lattice $Q\subset \ZZ^m$ with the sides $q_1,\dots,q_m$, respectively,
we say that $\{z_I\, :\, I\in Q\}$ satisfies $a$ if for any $J\in \ZZ^m$ the
minimum $\min_I\{a_I+
z_{I+J}\}$ is attained at least twice, provided that $I+J\in Q$ for each $I$ from the support of $a$. Clearly, one could consider an arbitrary subset of $\ZZ^m$ rather than just a parallelepiped.

Similar to section~\ref{five} define $d_{q_1,\dots,q_m}:=d_{q_1,\dots,q_m}(a)$ to be the
dimension of the tropical linear prevariety $\{z_I\, :\, I\in Q\}$ satisfying $a$ (i.~e. a set of points satisfying a finite family of tropical linear equations \cite{MS}).
Analogously one defines $m_{q_1,\dots,q_m}$ with respect to tropical multivariable minimal recurrent sequences. There exist limits
$$H(a):=\lim_{q_1,\dots,q_m \to \infty} d_{q_1,\dots,q_m}/(q_1\cdots q_m), \quad
h(a):=\lim_{q_1,\dots,q_m \to \infty} m_{q_1,\dots,q_m}/(q_1\cdots q_m)$$
\noindent which we also call the {\it tropical entropy} (respectively, the {\it tropical minimal entropy}) of $a$. Again $0\le h(a)\le H(a)$. One can prove the following proposition similarly to the proof of Proposition~\ref{upper}.

\begin{proposition}
Let the support of $a$ be located in a cube with the side $r$. Then $H(a)\le 1-1/r^m$.
\end{proposition}

\begin{example}
Let $a$ be a vector with $a_I=0$ for all $I=(i_1,\dots,i_m),\, 1\le i_1,\dots,i_m\le r$, so its support is a cube with the side $r-1$. Then the following $z=\{z_I\}$ is a tropical multivariable recurrent sequence satisfying $a$. Put $z_{l_1r,l_2r,\dots,l_mr}=z_{l_1r+1,l_2r,\dots,l_mr}=0$ for all integers $l_1,\dots,l_m$, and the rest of the coordinates of $z$ being arbitrary non-negative. Therefore, $H(a)\ge 1-2/r^m$.
\end{example}

\vspace{2mm}

{\bf Acknowledgements}. The author is grateful to the grant RSF 16-11-10075,
 to MCCME for inspiring atmosphere, to A.~Tiskin and to an anonymous referee for useful remarks.

\end{document}